\documentclass[12pt]{amsart}
\usepackage[margin=2.0cm]{geometry}

\usepackage{amsfonts}
\usepackage{mathrsfs}
\usepackage{cite}
\usepackage{graphicx}

\newcommand{\R}{{\mathbb R}}

\newtheorem{theorem}{Theorem}[section]
\newtheorem{lemma}[theorem]{Lemma}

\newtheorem{corollary}[theorem]{Corollary}

\title{The isodiametric problem on the sphere and in the hyperbolic space}
\author{K\'aroly J. B\"or\"oczky, \'Ad\'am Sagmeister }
\address{Alfr\'ed R\'enyi Institute of Mathematics, Hungarian Academy
  of Sciences, Reltanoda u. 13-15, H-1053 Budapest, Hungary, and
Department of Mathematics, Central European University, Nador u 9, H-1051, Budapest, Hungary} \email{boroczky.karoly.j@renyi.mta.hu}
\address{Lorant E\"otv\"os University, Institute of Mathematics, P\'azm\'any P\'eter s\'et\'any 1/c, Budapest, H-1117 Hungary} \email{sagmeister.adam@gmail.com }
\subjclass[2010]{Primary: }

\begin{document}

\maketitle

\section{Introduction}

Let $\mathcal{M}^n$ be either the Euclidean space $\R^n$, hyperbolic space $H^n$ or spherical space $S^n$ for $n\geq 2$. We write $V_{\mathcal{M}^n}$ to denote the $n$-dimensional volume (Lebesgue measure) on $\mathcal{M}^n$,
and $d_{\mathcal{M}^n}(x,y)$ to denote the geodesic distance between $x,y\in\mathcal{M}^n$.

For a bounded set $X\subset \mathcal{M}^n$,  its diameter ${\rm diam}_{\mathcal{M}^n} X$ is the supremum of the geodesic distances $d_{\mathcal{M}^n}(x,y)$ for $x,y\in X$.  
 For $D>0$ and $n\geq 1$, our goal is to determine the maximal volume of a subset of  $\mathcal{M}^n$ of diameter at most $D$. For any $z\in \mathcal{M}^n$ and $r>0$, let
$$
B_{\mathcal{M}^n}(z,r)=\{x\in \mathcal{M}^n:\, d(x,z)\leq r\}
$$ 
be the $n$-dimensional ball centered at $z$ where it is natural to assume $r<\pi$ if $\mathcal{M}^n=S^n$.
When it is clear from the context what space we consider, we drop the subscript referring to the ambient space.
In order to speak about the volume of a ball of radius $r$, we fix a reference point $z_0\in \mathcal{M}^n$ where
$z_0=o$ the origin if $\mathcal{M}^n=\R^n$.

It is well-known, due to Bieberbach \cite{Bie15} in $\R^2$ and 
P. Urysohn \cite{Ury24} in $\R^n$, that if $X\subset \R^n$ is measurable and bounded with 
${\rm diam}X=D$, then
\begin{equation}
\label{eucklidean}
V(X)\leq V(B(o,D/2)),
\end{equation}
and equality holds if and only if the closure of $X$ is a ball of radius $D/2$.

We prove the following hyperbolic- and spherical analogues of (\ref{eucklidean}).

\begin{theorem}
\label{hyperbolic}
If $D>0$ and $X\subset H^n$ is measurable and bounded with 
${\rm diam}X\leq D$, then
$$
V(X)\leq V(B(z_0,D/2)),
$$
and equality holds if and only if the closure of $X$ is a ball of radius $D/2$.
\end{theorem}

\begin{theorem}
\label{spherical}
If $D\in(0,\pi)$ and $X\subset S^n$ is measurable with 
${\rm diam}X\leq D$, then
$$
V(X)\leq V(B(z_0,D/2)),
$$
and equality holds if and only if the closure of $X$ is a ball of radius $D/2$.
\end{theorem}

On $S^2$, Hern\'andez Cifre, Mart\'{\i}nez Fern\'andez \cite{HCF15} proved a stronger version of
Theorem~\ref{spherical} for centrally symmetric sets of diameter less than $\pi/2$.

The proofs of Theorem~\ref{hyperbolic} and Theorem~\ref{spherical} build strongly on ideas related to two-point symmetrization in the paper
Aubrun, Fradelizi \cite{AuF04}. 
After reviewing the basic properties of spaces of constant curvature in Secion~\ref{secconstant}, 
we consider the extremal sets in Section~\ref{secDmaximal} and convex  sets in
Section~\ref{secconvex} from our point of view.  The main tool of this paper, two-point symmetrization, is introduced in
Section~\ref{sectwopoint} where we actually prove Theorem~\ref{spherical} if $D\leq\frac{\pi}2$ and Theorem~\ref{hyperbolic}. The reason why the argument is reasonably simple (kind of a proof from Erd\H{o}s' Book) for
 Theorem~\ref{spherical} if $D\leq\frac{\pi}2$ and for Theorem~\ref{hyperbolic} is that it is easy to show that the extremal sets are convex in these cases. However, if $D>\frac{\pi}2$ in the spherical case, then {\it a priori} much less information is availabe on the extremal sets, therefore a more technical argument is provided in
Section~\ref{secnonconvexproof}.

\section{Spaces of constant curvature}
\label{secconstant}

Let $\mathcal{M}^n$ be either $\R^n$,  $H^n$ or $S^n$. Our focus is on the spherical- and hyperbolic space, and we assume that $S^n$ is embedded into $\R^{n+1}$ the standard way,
and $H^n$ is embedded into $\R^{n+1}$ using the hyperboloid model. We write $\langle \cdot,\cdot\rangle$ to denote
the standard scalar product in $\R^{n+1}$, write 
$z^\bot=\{x\in\R^{n+1}:\,\langle x,e\rangle=0\}$ for $z\in\R^{n+1}\backslash o$ and fix an $e\in\R^{n+1}$.
In particular, we have
\begin{eqnarray*}
S^n&=& \{x\in\R^{n+1}:\,\langle x,x\rangle=1\}\\
H^n&=& \{x+te:\,x\in e^\bot \mbox{ and }t\geq 1\mbox{ and }t^2-\langle x,x\rangle=1\}.
\end{eqnarray*}

For $H^n$, we also consider the following symmetric bilinear form $\mathcal{B}$ on $\R^{n+1}$: If $x=x_0+te\in\R^{n+1}$ and $y=y_0+se\in\R^{n+1}$ for
$x_0,y_0\in e^\bot$ and $t,s\in\R$, then
$$
\mathcal{B}(x,y)=ts-\langle x_0,y_0\rangle.
$$
In particular,
\begin{equation}
\label{Hn}
\mathcal{B}(x,x)=1\mbox{ for $x\in H^n$}.
\end{equation}

Again let $\mathcal{M}^n$ be either $\R^n$,  $H^n$ or $S^n$ using the models as above for $H^n$ and $S^n$. For
$z\in \mathcal{M}^n$, we define the tangent space $T_z$ as
\begin{eqnarray*}
T_z&=&\{x\in\R^{n+1}:\,\mathcal{B}(x,z)=0\} \mbox{ if $\mathcal{M}^n=H^n$}\\
T_z&=&z^\bot \mbox{ if $\mathcal{M}^n=S^n$}\\
T_z&=&\R^n \mbox{ if $\mathcal{M}^n=\R^n$}.
\end{eqnarray*}  
We observe that $T_z$ is an $n$-dimensional real vector space equipped with the scalar product
$-\mathcal{B}(\cdot,\cdot)$ if $\mathcal{M}^n=H^n$, and with the scalar product 
$\langle \cdot,\cdot\rangle$ if $\mathcal{M}^n=S^n$ or $\mathcal{M}^n=\R^n$.

Let us consider lines and $(n-1)$-dimensional (totally geodesic) subspaces of $\mathcal{M}^n$. A line $\ell\subset \mathcal{M}^n$ passing through 
a $z\in \mathcal{M}^n$ is given by a unit vector $u\in T_z$; namely, $-\mathcal{B}(u,u)=1$ if 
$\mathcal{M}^n=H^n$, and $\langle u,u\rangle=1$ if 
$\mathcal{M}^n=S^n$ or $\mathcal{M}^n=\R^n$, and the line $\ell$ is parameterized by
\begin{eqnarray*}
x=({\rm ch}\,t)z+({\rm sh}\,t)u &\mbox{ if }&\mathcal{M}^n=H^n \mbox{ and $t\in \R$ }\\
x=(\cos t)z+(\sin t)u &\mbox{ if }&\mathcal{M}^n=S^n \mbox{ and $t\in [-\pi,\pi]$ }\\
x=z+tu &\mbox{ if }&\mathcal{M}^n=\R^n \mbox{ and $t\in \R$ }
\end{eqnarray*}
where $d_{\mathcal{M}^n}(x,z)=|t|$. In addition, the $(n-1)$-dimensional subspace of $\mathcal{M}^n$
passing through $z\in \mathcal{M}^n$ and having normal vector $v\in T_z\backslash \{o\}$ is
\begin{eqnarray*}
\{x\in H^n:\,\mathcal{B}(x,v)= 0\} &\mbox{ if }&\mathcal{M}^n=H^n \\
\{x\in S^n:\,\langle x,v\rangle = 0\} &\mbox{ if }&\mathcal{M}^n=S^n\\
\{x\in \R^n:\,\langle x,v\rangle = 0\} &\mbox{ if }&\mathcal{M}^n=\R^n.  
\end{eqnarray*}

We say that a non-empty compact set $C\subset \mathcal{M}^n$ with $C\neq S^n$ is solid, if $C$ is the closure of ${\rm int}_{\mathcal{M}^n} C$.
For such a $C$, we say that $x\in\partial_{\mathcal{M}^n}C$ is strongly regular if there exists
$r>0$ and $z\in C$ such that $B_{\mathcal{M}^n}(z,r)\subset C$ and 
$x\in\partial_{\mathcal{M}^n} B_{\mathcal{M}^n}(z,r)$. In this case, 
we set $N_C(x)\in T_x$ to be an exterior unit vector at $x$; namely,
\begin{eqnarray*}
-\mathcal{B}(N_C(x),N_C(x))&=&1 \mbox{ and $z=({\rm ch}\,r)x-({\rm sh}\,r)N_C(x)$ if $\mathcal{M}^n=H^n$}\\
\langle N_C(x),N_C(x)\rangle &=&1 \mbox{ and $z=(\cos r)x-(\sin r)N_C(x)$ if $\mathcal{M}^n=S^n$}\\
\langle N_C(x),N_C(x)\rangle &=&1 \mbox{ and $z=x- rN_C(x)$ if $\mathcal{M}^n=\R^n$}.
\end{eqnarray*} 

\begin{lemma}
\label{strongly regular-dense}
If $C\subset \mathcal{M}^n$ is solid, then the strongly regular boundary points are dense 
in $\partial_{\mathcal{M}^n}C$.
\end{lemma}
\proof
Let $z\in\partial_{\mathcal{M}^n}C$, and let $\varepsilon>0$. As $C$ is solid, we may choose
an $x\in{\rm int}_{\mathcal{M}^n}C$ such that $d_{\mathcal{M}^n}(x,z)<\varepsilon/2$. Let $r>0$
be maximal with the property that $B_{\mathcal{M}^n}(x,r)\subset C$, and hence
$r\leq d_{\mathcal{M}^n}(x,z)$ and there exists
$$
y\in \partial_{\mathcal{M}^n}B_{\mathcal{M}^n}(z,r)\cap \partial_{\mathcal{M}^n} C.
$$
Therefore $y$ is a strongly regular boundary point, and the triangle inequality yields that $d_{\mathcal{M}^n}(y,z)<\varepsilon$.
\endproof

According to Vinberg \cite{Vin93}, the ``Standard Hypersurfaces" in either $\R^n$, $H^n$ or $S^n$ are as follows:
\begin{itemize}
\item Hyperplanes of the form
$\{x\in \R^n:\,\langle x,p\rangle=t\}$ in $\R^n$ for $p\in \R^n\backslash o$ and $t\in\R$;
\item $\partial_{\R^n} B(z,r)=\{x\in \R^n:\,\langle x-z,x-z\rangle=r^2\}$ in $\R^n$ for $z\in \R^n$ and $r>0$;
\item $\partial_{S^n} B(z,r)=\{x\in S^n:\,\langle x,z\rangle=\cos r\}$ in $S^n$ for $z\in S^n$ and $r\in(0,\pi)$,  hence
for any $p\in \R^{n+1}\backslash o$ and $t\in \R$, the set $\{x\in S^n:\,\langle x,p\rangle=t\}$ is either empty, a point, or the boundary of a spherical ball;
\item $\partial_{H^n} B(z,r)=\{x\in H^n:\,\langle x,z\rangle={\rm ch}\, r\}$ in $H^n$ for 
$r>0$ and $z\in H^n$, hence  for any $t\in \R$ and $p\in \R^{n+1}\backslash o$
with $\mathcal{B}(p,p)>0$, the set $\{x\in H^n:\,\mathcal{B}(x,p)=t\}$ is either empty, a point, or the boundary of a hyperbolic ball;
\item Hyperplanes of the form
$\{x\in H^n:\,\mathcal{B}(x,p)=0\}$ in $H^n$ for $p\in \R^n\backslash o$ with $\mathcal{B}(p,p)<0$;
\item Hyperspheres of the form
$\{x\in H^n:\,\mathcal{B}(x,p)=t\}$ in $H^n$ for $t\in \R\backslash 0$ and $p\in \R^n\backslash o$ with $\mathcal{B}(p,p)<0$;
\item Horospheres of the form
$\{x\in H^n:\,\mathcal{B}(x,p)=t\}$ in $H^n$ for $t\in \R\backslash 0$ and $p\in \R^n\backslash o$ with $\mathcal{B}(p,p)=0$.
\end{itemize}

We note the following properties.

\begin{lemma}
\label{hypersurf-component}
If $\Xi$ is a standard hypersurface in $\mathcal{M}^n$
where $\mathcal{M}^n$ is either $\R^n$, $H^n$ or $S^n$, then
$\mathcal{M}^n\backslash\Xi$ has two connected components,
the boundary of both components is $\Xi$, and any of these components is bounded if and only if the component is an open ball.
\end{lemma} 

\begin{corollary}
\label{hypersurf-solid}
If $\mathcal{M}^n$ is either $\R^n$, $H^n$ or $S^n$, and $C\subset \mathcal{M}^n$, $C\neq\mathcal{M}^n$ is a solid set whose strongly regular boundary points are contained in a fixed standard hypersurface, then $C$ is a ball.
\end{corollary}
\proof Let $\Xi$ be  standard hypersurface containing the strongly regular boundary points of $C$. 
Since $\Xi$ is closed and strongly regular boundary points are dense in $\partial_{\mathcal{M}^n} C$ according to Lemma~\ref{strongly regular-dense}, we have
\begin{equation}
\label{inXi}
\partial_{\mathcal{M}^n} C\subset \Xi.
\end{equation}
Next we prove
\begin{equation}
\label{equalXi}
\partial_{\mathcal{M}^n} C= \Xi.
\end{equation}
We suppose that there exists $z\in \Xi\backslash \partial_{\mathcal{M}^n} C$, and seek a contradiction. We consider
some $y\in({\rm int}_{\mathcal{M}^n}C)\backslash \Xi$, thus Lemma~\ref{hypersurf-component} yields that
there exists a continuous curve $\gamma:[0,1]\to \mathcal{M}^n$ such that $\gamma(0)=y$, $\gamma(1)=z$
and $\gamma(t)\not\in \Xi$ for $t<1$. As $z\not\in C$, there exists $s=\max\{t:\,\gamma(t)\in C\}<1$. It follows that
$\gamma(s)\in (\partial_{\mathcal{M}^n} C)\backslash \Xi$, contradicting (\ref{inXi}), and proving (\ref{equalXi}).

We deduce from (\ref{equalXi}) that
${\rm int}_{\mathcal{M}^n} C$ is the union of the components of $\mathcal{M}^n\backslash\Xi$,
and since $C\neq\mathcal{M}^n$, 
${\rm int}_{\mathcal{M}^n} C$ is one of the components of $\mathcal{M}^n\backslash\Xi$ by Lemma~\ref{hypersurf-component}.
As $C$ is bounded, we conclude that $C$ is a ball again by Lemma~\ref{hypersurf-component}.
\endproof

In the final part of this section, we assume that $\mathcal{M}^n$ is either  $H^n$ or $S^n$, and 
use the models in $\R^{n+1}$ above.
For $k=1,\ldots,n-1$, the $k$-dimensional totally geodesic subspaces are of the form $L\cap \mathcal{M}^n$ where
$L$ is a linear $(k+1)$-dimensional subspace of $\R^{n+1}$ with $L\cap \mathcal{M}^n\neq\emptyset$. Next, we define $\pi:\R^{n+1}\backslash e^\bot\to e^\bot+e$ by
$$
\pi(x)=\frac{x}{\langle x,e\rangle}.
$$
It follows that the restriction of $\pi$ to $H^n$ is a diffeomorphism into
the ``open" $n$-ball $\{x\in e^\bot+e:\,d_{\R^{n+1}}(x,e)<1\}$, and 
the restriction of $\pi$ to ${\rm int}_{S^n}B_{S^n}(e,\frac{\pi}2)$ is a diffeomorphism into
the affine $n$-plane $e^\bot+e$ of $\R^{n+1}$. In addition, 
for any $k=1,\ldots,n-1$, $\pi$ induces a natural bijection between
certain $k$-dimensional affine subspaces of $e^\bot+e$ and certain 
$k$-dimensional totally geodesic subspaces of $\mathcal{M}^n$ not contained in $e^\bot$.
In particular, if $L$ is a $(k+1)$-dimensional linear subspace of $\R^{n+1}$, then
\begin{eqnarray}
\label{piHn}
\pi(L\cap H^n)&=&L\cap(e^\bot+e) \mbox{ \ provided $L\cap H^n\neq \emptyset$}\\
\label{piSn}
\pi\left(L\cap {\rm int}_{S^n}B_{S^n}\left(e,\frac{\pi}2\right)\right)&=&L\cap(e^\bot+e) \mbox{ \ provided 
$L\not\subset e^\bot$.}
\end{eqnarray}

\section{$D$-maximal sets}
\label{secDmaximal}

The main tool to obtain convex bodies with extremal properties is the Blaschke Selection Theorem. 
First we impose a metric on non-empty compact subsets. Let $\mathcal{M}^n$ be either $\R^n$,  $H^n$ or $S^n$.
For a non-empty compact set $C\subset\mathcal{M}^n$ and $z\in \mathcal{M}^n$, we set $d_{\mathcal{M}^n}(z,C)=\min_{x\in C}d_{\mathcal{M}^n}(z,x)$.
For any non-empty compact set $C_1,C_2\subset\mathcal{M}^n$, we define their Hausdorff distance
$$
\delta_{\mathcal{M}^n}(C_1,C_2)=\max\left\{
\max_{x\in C_2}d_{\mathcal{M}^n}(x,C_1),
\max_{y\in C_1}d_{\mathcal{M}^n}(y,C_2)\right\}.
$$
The Hausdorff distance is a metric on the space of non-empty compact subsets in $\mathcal{M}^n$.
We say that a sequence $\{C_m\}$ of non-empty compact subsets of $\mathcal{M}^n$ is bounded if there is a ball containing every $C_m$.
For non-empty compact sets $C_m,C\subset \mathcal{M}^n$, we write $C_m\to C$ to denote if
the sequence $\{C_m\}$ tends to $C$ in terms of the Hausdorff distance. 

The following is well-known and we present the easy argument for the sake of completeness (see Theorem~1.8.8 in R. Schneider \cite{Sch14} for the case when each set is convex).

\begin{lemma}
\label{Hausdorffconvergence}
For non-empty compact sets $C_m,C\subset \mathcal{M}^n$ where $\mathcal{M}^n$ is either $\R^n$,  $H^n$ or $S^n$, we have $C_m\to C$ if and only if
\begin{description}
\item[(i)] assuming $x_m\in C_m$, the sequence $\{x_m\}$ is bounded and any accumulation point of  $\{x_m\}$  
lies in $C$;
\item[(ii)] for any $y\in C$, there exist $x_m\in C_m$ for each $m$ such that $\lim_{m\to \infty}x_m=y$.
\end{description}
\end{lemma}
\proof First we assume that $C_m\to C$. For (i), if $x_m\in C_m$, then let $y_m\in C$ be a point closest to $x_m$. Now $\lim_{m\to\infty}d_{\mathcal{M}^n}(x_m,C_m)=0$ yields that $\{x_m\}$ is bounded and the sequences $\{x_m\}$ and $\{y_m\}$ have the same set of accumulation points, proving (i). For (ii), let $y\in C$, and let $x_m\in C_m$ be a point closest to $y$.
We have $\lim_{m\to \infty}x_m=y$ because $\lim_{m\to \infty}d_{\mathcal{M}^n}(y,C_m)=0$. 

Next we assume that (i) and (ii) hold for $\{C_m\}$ and $C$. Let  $x_m\in C_m$ be a point farthest from $C$ and let
$y_m\in C$ be a point farthest to $C_m$, and hence both $\{x_m\}$ and $\{y_m\}$ are bounded. We choose
subsequences $\{x_{m'}\}\subset\{x_m\}$ and $\{y_{m"}\}\subset\{y_m\}$ such that
$$
\lim_{m'\to\infty} d_{\mathcal{M}^n}(x_{m'},C)=\limsup_{m\to\infty} d_{\mathcal{M}^n}(x_m,C)
\mbox{ \ and \ }
\lim_{m"\to\infty} d_{\mathcal{M}^n}(y_{m"},C_{m"})=\limsup_{m\to\infty} d_{\mathcal{M}^n}(y_m,C_m).
$$
We may also assume that $\lim_{m'\to\infty}x_{m'}=x$ and $\lim_{m"\to\infty}y_{m"}=y$ where
$x\in C$ by (i) and $y\in C$ by the compactness of $C$. On the one hand, it follows that
$\lim_{m\to\infty} d_{\mathcal{M}^n}(x_m,C)=0$. On the other hand, there exists $z_m\in C_m$ such that
$\lim_{m\to\infty}z_m=y$ by (ii), therefore $\lim_{m\to\infty} d_{\mathcal{M}^n}(y_m,C_m)=0$ as well.
\endproof

The space of non-empty compact subsets of $\mathcal{M}^n$ is locally compact according to the Blaschke Selection Theorem
(see R. Schneider \cite{Sch14}).

\begin{theorem}[Blaschke]
\label{Blaschke}
If $\mathcal{M}^n$ is either $\R^n$,  $H^n$ or $S^n$,
then any bounded sequence of non-empty compact subsets of $\mathcal{M}^n$ has
a convergent subsequence.
\end{theorem}

Let us consider convergent sequences of compact subsets of $\mathcal{M}^n$.

\begin{lemma}
\label{limit}
Let $\mathcal{M}^n$ be either $\R^n$,  $H^n$ or $S^n$, and
let the sequence $\{C_m\}$ of non-empty compact subsets of $\mathcal{M}^n$ tend to $C$.
\begin{description}
\item[(i)] ${\rm diam}_{\mathcal{M}^n}\,C=\lim_{m\to \infty}{\rm diam}_{\mathcal{M}^n}\,C_m$
\item[(ii)] $V_{\mathcal{M}^n}(C)\geq \limsup_{m\to \infty} V_{\mathcal{M}^n}(C_m)$
\end{description}
\end{lemma}
\proof (i) follows from Lemma~\ref{Hausdorffconvergence}.

For (ii), it is sufficient to prove that for any $\varepsilon>0$, there exists $M$ such that
$V_{\mathcal{M}^n}(C_m)\leq V_{\mathcal{M}^n}(C)+\varepsilon$ if $m\geq M$. Choose $r>0$
such that the open set
$$
U_r=\bigcup_{x\in C}{\rm int}_{\mathcal{M}^n}B_{\mathcal{M}^n}(x,r)
$$
satisfies $V_{\mathcal{M}^n}(U_r)\leq V_{\mathcal{M}^n}(C)+\varepsilon$. Such an $r$ exists as $C$ is compact.
We choose $M$ such that $\delta_{\mathcal{M}^n}(C_m,C)<r$ if $m\geq M$, and hence
$C_m\subset U_r$ if $m\geq M$.
\endproof

For $D>0$ where $D<\pi$ if $\mathcal{M}^n=S^n$, we say that a compact set $C\subset \mathcal{M}^n$ 
is $D$-maximal if ${\rm diam}_{\mathcal{M}^n}\,C\leq D$ and
$$
V_{\mathcal{M}^n}(C)=\sup\{V_{\mathcal{M}^n}(X):\, X\subset \mathcal{M}^n\mbox{ compact and }
{\rm diam}_{\mathcal{M}^n} X\leq D \}.
$$

\begin{theorem}
\label{maximal-sets}
Let $\mathcal{M}^n$ be either $\R^n$,  $H^n$ or $S^n$, and let
 $D>0$ where $D<\pi$ if $\mathcal{M}^n=S^n$.
\begin{description}
\item[(i)] There exists a $D$-maximal set in $\mathcal{M}^n$.
\item[(ii)] For any $D$-maximal set $C$ in $\mathcal{M}^n$ and $z\in\partial_{\mathcal{M}^n}C$,
there exists $y\in\partial_{\mathcal{M}^n}C$ such that $d_{\mathcal{M}^n}(z,y)=D$.
\end{description}
\end{theorem}
\proof Let $\{C_m\}$ be a sequence of compact subsets of $\mathcal{M}^n$ with $z_0\in C_m$,
${\rm diam}_{\mathcal{M}^n}\,C_m\leq D$ and
$$
\lim_{m\to \infty}V_{\mathcal{M}^n}(C_m)
=\sup\{V_{\mathcal{M}^n}(X):\, X\subset \mathcal{M}^n\mbox{ compact and }
{\rm diam}_{\mathcal{M}^n} X\leq D \}.
$$
According to the Blaschke Selection Theorem Theorem~\ref{Blaschke}, we may asume that
the sequence $\{C_m\}$ tends to a  compact subset $X\subset \mathcal{M}^n$. Here $X$ is a $D$-maximal set by Lemma~\ref{limit}.

Next let $C$ be any $D$-maximal set in $\mathcal{M}^n$, and let $z\in\partial_{\mathcal{M}^n}C$.
We suppose that \\
$\Delta=\max_{x\in C}d_{\mathcal{M}^n}(z,x)<D$,
and seek a contradiction. As $z\in\partial_{\mathcal{M}^n}C$, there exists some 
$y\in B_{\mathcal{M}^n}(z,\frac12(D-\Delta))\backslash C$, and hence
$B_{\mathcal{M}^n}(y,r)\cap C=\emptyset$ for some $r\in(0,\frac12(D-\Delta))$.
Therefore $X=C\cup  B_{\mathcal{M}^n}(z,r)$ satisfies that
$V_{\mathcal{M}^n}(X)>V_{\mathcal{M}^n}(C)$ and ${\rm diam}_{\mathcal{M}^n}\,X\leq D$,
which is a contradiction verifying (ii).
\endproof

\section{Convex sets}
\label{secconvex}

Let $\mathcal{M}^n$ be either $\R^n$,  $H^n$ or $S^n$.
For $x,y\in \mathcal{M}^n$ where $x\neq - y$ if $\mathcal{M}^n=S^n$, we write $[x,y]_{\mathcal{M}^n}$ to denote the geodesic segment between $x$ and $y$ whose length is 
$d_{\mathcal{M}^n}(x,y)$. 
We call $X\subset \mathcal{M}^n$ convex if $[x,y]_{\mathcal{M}^n}\subset X$
for any $x,y\in X$, and in addition we assume that $X$ is contained in an open hemisphere if $\mathcal{M}^n=S^n$.
 For $Z\subset \mathcal{M}^n$ where we assume that 
$Z$ is contained in an open hemisphere if $\mathcal{M}^n=S^n$, the convex hull ${\rm conv}_{\mathcal{M}^n} Z$ is the intersection of all convex sets containing $Z$. 

We observe that for a non-empty compact convex $Z\subset\mathcal{M}^n$, the conditions that $V_{\mathcal{M}^n}(Z)>0$,
${\rm int}_{\mathcal{M}^n}Z\neq\emptyset$ and $Z$ is solid are equivalent.

We deduce from (\ref{piHn}) and (\ref{piSn}) that
if $x,y\in H^n$ or   $x,y\in {\rm int}_{S^n}B_{S^n}(e,\frac{\pi}2)$,
and $\mathcal{M}^n$ is either $H^n$ or $S^n$, respectively, then 
$\pi([x,y]_{\mathcal{M}^n})=[\pi(x),\pi(y)]_{\R^{n+1}}$; namely, the Euclidean segment in $e^\bot+e$.
Thus for a 
subset $Z$  of either $H^n$ or ${\rm int}_{S^n}B_{S^n}(e,\frac{\pi}2)$, $Z$ is convex
on $H^n$ or $S^n$, respectively, if and only if $\pi(Z)\subset e^\bot+e$ is convex.

Since on the sphere, it is an important issue whether a set $Z\subset S^n$ is contained in an open hemisphere, we note the following condition:

\begin{lemma}
\label{diameter-on-sphere}
If ${\rm diam}_{S^n}X<\arccos \frac{-1}{n+1}$ for $X\subset S^n$, then $X$ is contained in
an open hemisphere. 
\end{lemma}
\proof
We may assume that $X$ is compact. Let $Z={\rm conv}_{\R^{n+1}}\,X$, and hence $Z$ is compact as well. 
Let $z\in Z$ be the unique closest point of $Z$ to $o$.
It follows from the Charath\'eodory theorem applied in $\R^{n+1}$
that there exist $x_1,\ldots,x_{n+2}\in X$ and $\alpha_1,\ldots,\alpha_{n+2}\in[0,1]$
satisfying $\alpha_1+\ldots +\alpha_{n+2}=1$ and
$$
\alpha_1x_1+\ldots +\alpha_{n+2}x_{n+2}=z.
$$
As ${\rm diam}_{S^n}X<\arccos \frac{-1}{n+1}$, we have
$\langle x_i,x_j\rangle>\frac{-1}{n+1}$ for any $i\neq j$. We deduce
from $2\alpha_i\alpha_j\leq \alpha_i^2+\alpha_j^2$ that
$$
\langle z,z\rangle> \left(\sum_{i=1}^{n+2}\alpha_i^2\right)-
\left(\sum_{i< j}\frac{2\alpha_i\alpha_j}{n+1}\right)
\geq \left(\sum_{i=1}^{n+2}\alpha_i^2\right)-
\left(\sum_{i< j}\frac{\alpha_i^2+\alpha_j^2}{n+1}\right)=0
$$
therefore $X\subset Z\subset \{x\in R^{n+1}:\,\langle x,z\rangle>0\}$.
\endproof

\begin{lemma}
\label{ball-convex}
If either $\mathcal{M}^n=H^n$ and $r>0$, or $\mathcal{M}^n=S^n$
and $r\in(0,\frac{\pi}2)$, then $B_{\mathcal{M}^n}(z,r)$ is convex for any $z\in \mathcal{M}^n$. In addition,
if $X\subset S^n$ convex, then $X\cap B_{S^n}(z,\frac{\pi}2)$ is convex.
\end{lemma}
\proof For the case $\mathcal{M}^n=H^n$ and $r>0$, or $\mathcal{M}^n=S^n$
and $r\in(0,\frac{\pi}2)$, we may assume that $z=e$. Thus $\pi(B_{\mathcal{M}^n}(z,r))$ is a Euclidean ball in $e^\bot+e$, which in turn yields that $B_{\mathcal{M}^n}(e,r)$ is convex.

If $X\subset S^n$ is convex, then we may assume that $X\subset {\rm int}_{S^n}B_{S^n}(e,\frac{\pi}2)$.
For $H^+=\{x\in \R^{n+1}:\,\langle z,x\rangle\geq 0\}$, we have
$$
\pi\left(X\cap B_{S^n}\left(z,\frac{\pi}2\right)\right)=\pi(X)\cap H^+,
$$
which is convex, and hence $X\cap B_{S^n}(z,\frac{\pi}2)$ is convex as well.
\endproof

We remark that 
 $B_{S^n}(z,r)$ is not convex if $r\in[\frac{\pi}2,\pi)$.

\begin{lemma} 
\label{convexhull}
Let $\mathcal{M}^n$ be either $\R^n$,  $H^n$ or $S^n$, and let
$X\subset \mathcal{M}^n$ be compact, non-empty and satisfy
${\rm diam}\,X\leq \frac{\pi}2$ in the case of  $\mathcal{M}^n=S^n$. Then
\begin{description}
\item[(i)] ${\rm diam}_{\mathcal{M}^n}{\rm conv}_{\mathcal{M}^n}\,X= {\rm diam}_{\mathcal{M}^n}\,X$;
\item[(ii)] $V_{\mathcal{M}^n}\left({\rm conv}_{\mathcal{M}^n}\,X\right)>
V_{\mathcal{M}^n}(X)$ if $V_{\mathcal{M}^n}\left({\rm conv}_{\mathcal{M}^n}\,X\right)>0$
and ${\rm conv}_{\mathcal{M}^n}\,X\neq X$.
\end{description}
\end{lemma}
\proof 
For (i), let ${\rm diam}\,X=D$ and let $x_1,x_2\in {\rm conv}_{\mathcal{M}^n} X$.
First, we consider the case when $\mathcal{M}^n$ is either $\R^n$ or  $H^n$.
Since $X\subset B_{\mathcal{M}^n}(x_1,D)$ and  $B_{\mathcal{M}^n}(x_1,D)$ is convex
 by Lemma~\ref{ball-convex}, we have
$x_2\in B_{\mathcal{M}^n}(x_1,D)$. Therefore $d_{\mathcal{M}^n}(x_1,x_2)\leq D$.

If $\mathcal{M}^n=S^n$, then Lemma~\ref{diameter-on-sphere} yields that
$X\subset B_{S^n}(z,R)$ for some $z\in S^n$ and $R\in(0,\frac{\pi}2)$.
Since $X\subset B_{S^n}(x_1,D)\cap B_{S^n}(z,R)$, which is convex
by Lemma~\ref{ball-convex}, we have
$x_2\in B_{S^n}(x_1,D)\cap B_{S^n}(z,R)$. Therefore $d_{S^n}(x_1,x_2)\leq D$,
finally proving (i).

For (ii), we assume that $V(Z)>0$ for $Z={\rm conv}_{\mathcal{M}^n}\,X$ and $Z\neq X$.
As $Z$ is convex, it follows that the closure of ${\rm int} Z$ is $Z$
Since  $X$ is compact and $X\neq Z$, there exists some $z\in ({\rm int} Z)\backslash X$. Therefore
$B_{\mathcal{M}^n}(z,r)\subset  ({\rm int} Z)\backslash X$ for some $r>0$, proving that
$V(Z)>V(X)$.
\endproof

Theorem~\ref{maximal-sets} guarantees the existence of $D$-maximal sets, and Lemma~\ref{convexhull}
yields

\begin{corollary}
\label{maximal-convex}
If $\mathcal{M}^n$ is either $\R^n$,  $H^n$ or $S^n$, and 
 $D>0$ where $D\leq\frac{\pi}2$ if $\mathcal{M}^n=S^n$, then
 any $D$-maximal set in $\mathcal{M}^n$ is convex.
\end{corollary}

Using the map $\pi$, (\ref{piHn}) and (\ref{piSn}) in the spherical- and hyperbolic case, we deduce 
from Lemma~\ref{strongly regular-dense}

\begin{lemma}
\label{supporting-hyperplane}
If $\mathcal{M}^n$ is either $\R^n$,  $H^n$ or $S^n$, 
$K\subset \mathcal{M}^n$ is convex with ${\rm int}_{\mathcal{M}^n}K\neq\emptyset$
and $z\in \partial_{\mathcal{M}^n}K$, then there exists a supporting hyperplane of $\mathcal{M}^n$
containing $z$ and not intersecting ${\rm int}_{\mathcal{M}^n}K$. In addition, 
if $z$ is a strongly regular boundary point, then there exists a unique supporting hyperplane, and its
unit exterior normal is $N_K(z)$.
\end{lemma}
 
In turn, we deduce the following statement, which will be important in identifying boundary points of
a two-point symmetrization.

\begin{lemma}
\label{secant}
Let $\mathcal{M}^n$ be either $\R^n$,  $H^n$ or $S^n$, 
let $K\subset \mathcal{M}^n$ be convex with ${\rm int}_{\mathcal{M}^n}K\neq\emptyset$,
and let $x,y\in \partial_{\mathcal{M}^n}K$.
\begin{description}
\item[(i)] If $[x,y]_{\mathcal{M}^n}\cap{\rm int}_{\mathcal{M}^n}K\neq\emptyset$,
then $\ell\cap K=[x,y]_{\mathcal{M}^n}$ for the line $\ell$ passing through $x$ and $y$.
\item[(ii)] If $[x,y]_{\mathcal{M}^n}\subset\partial_{\mathcal{M}^n}K$, then there exists
a supporting hyperplane $\Pi$ to $K$ containing $[x,y]_{\mathcal{M}^n}$.
\end{description}
\end{lemma}

\section{Two-point symmetrization}
\label{sectwopoint}

Let $\mathcal{M}^n$ be either $\R^n$,  $H^n$ or $S^n$, let $H^+$ be a closed half space bounded by the $(n-1)$-dimensional subspace $H$ in $\mathcal{M}^n$, and let
$X\subset \mathcal{M}^n$ be compact and non-empty. We write $H^-$ to denote the other closed half space of $\mathcal{M}^n$ determined by $H$ and
$\sigma_H(X)$ to denote the reflected image of $X$ through the $(n-1)$-subspace $H$. 

The two-point symmetrization $\tau_{H^+}(X)$ of $X$ with respect to $H^+$ is a rearrangement of $X$ by replacing
$(H^-\cap X)\backslash \sigma_H(X)$ by its reflected image through $H$ where readily this reflected image is disjoint from $X$. Naturally, interchanging the role of $H^+$ and $H^-$ results in taking the reflected image of $\tau_{H^+}(X)$ through $H$. Since this operation does not change any relevant property of the new set, we simply use the notation  $\tau_H(X)$ (see Figure 1).

\begin{figure}
\begin{center}
\includegraphics[width=18em]{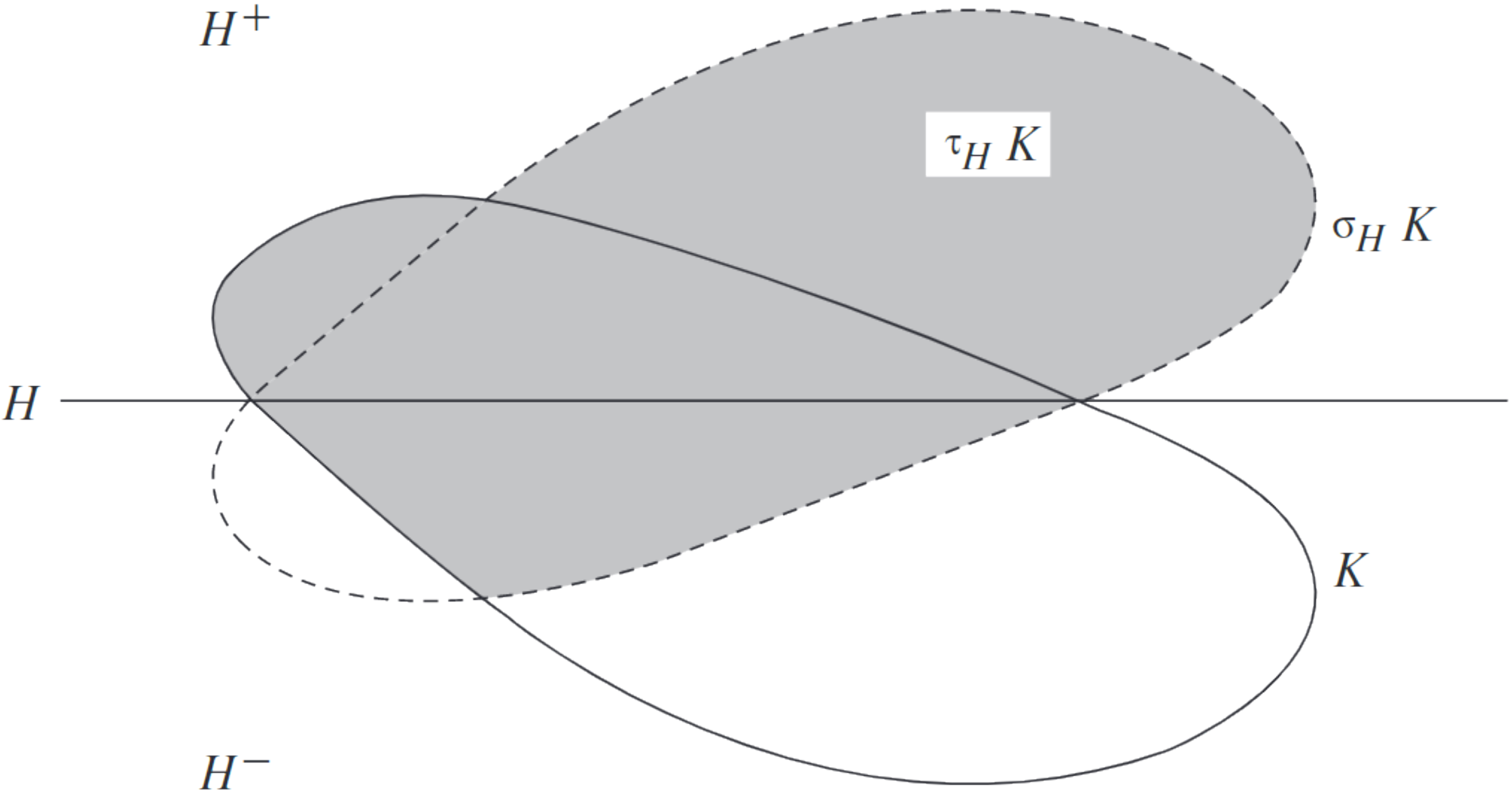}
\end{center}
\caption{ }
\end{figure}  

\begin{lemma}
\label{two-point-prop}
Let $\mathcal{M}^n$ be either $\R^n$,  $H^n$ or $S^n$, let $H^+$ be a half space, and let
$X\subset \mathcal{M}^n$ be compact and non-empty such that ${\rm diam}_{\mathcal{M}^n}(X)<\pi$ if $\mathcal{M}^n=S^n$. Then $\tau_H(X)=\tau_{H^+}(X)$ is compact and satisfies
\begin{description}
\item[(i)] $\tau_H(X)\cap H^+=\big(X\cup\sigma_H(X)\big)\cap H^+$;
\item[(ii)] $\tau_H(X)\cap H^-=\big(X\cap\sigma_H(X)\big)\cap H^-$;
\item[(iii)] $V_{\mathcal{M}^n}(\tau_H(X))=V_{\mathcal{M}^n}(X)$;
\item[(iv)] ${\rm diam}_{\mathcal{M}^n}(\tau_H(X))\leq {\rm diam}_{\mathcal{M}^n}(X)$.
\end{description}
\end{lemma}
\proof Properties (i) and (ii) are just reformulations of the definition of two-point symmetrization, and they
directly yield (iii) and the compactness of $\tau_H(X)$.

For (iv), let $x,y\in \tau_H(X)$.
 If either $x,y\in X$ or $x,y\in\sigma_H(X)$, then readily 
$d_{\mathcal{M}^n}(x,y)\leq {\rm diam}_{\mathcal{M}^n}(X)$. Otherwise, we may assume that
$x\in X\backslash \sigma_H(X)$ and $y\in  \sigma_H(X)\backslash X$, thus
$x\in X\cap (H^+\backslash  H)$ and $y\in \sigma_H(X)\cap (H^+\backslash H)$. 
It follows that $\sigma_H(y)\in X\cap (H^-\backslash H)$, and hence
$[x,\sigma_H(y)]_{\mathcal{M}^n}$ intersects $H$ in a unique $z$ where
$x,\sigma_H(y)\in X$ implies that $[x,\sigma_H(y)]_{\mathcal{M}^n}$ is well-defined even if
$\mathcal{M}^n=S^n$. 
Applying the triangle inequality to $x,y,z$, we deduce  that
\begin{eqnarray*}
d_{\mathcal{M}^n}(x,y)&\leq& d_{\mathcal{M}^n}(x,z)+d_{\mathcal{M}^n}(z,y)
=d_{\mathcal{M}^n}(x,z)+d_{\mathcal{M}^n}(z,\sigma_H(y))\\
&=&
d_{\mathcal{M}^n}(x,\sigma_H(y))\leq {\rm diam}_{\mathcal{M}^n}(X)
.
\end{eqnarray*}
\endproof

Two-point symmetrization appeared first in V. Wolontis \cite{Wol52}. It is applied to prove the
isoperimetric inequality in the spherical space by 
Y. Benyamini \cite{Ben84}, and the spherical analogue of the Blaschke-Santal\'o inequality by
Gao, Hug, Schneider \cite{GHS03}.

The following statement was proved by 
Aubrun, Fradelizi \cite{AuF04} in the Euclidean- and spherical case. Here we provide a somewhat more detailed version of the argument in \cite{AuF04}.
Part of the reason for more details is that 
in the case $\frac{\pi}2<D<\pi$ of Theorem~\ref{spherical} concerned about $S^n$,
 we use many ideas of the proof of Theorem~\ref{two-point-convex} on the one hand,  
 however, many essential ingredients of the argument for Theorem~\ref{two-point-convex} do not hold anymore.

\begin{theorem}
\label{two-point-convex}
If $\mathcal{M}^n$ is either $H^n$ or $S^n$, and the compact, convex $K\subset \mathcal{M}^n$
with non-empty interior satisfies that
$\tau_H(K)$ is convex for any $(n-1)$-dimensional subspace $H$ of $\mathcal{M}^n$, then $K$ is a ball.
\end{theorem}
\proof Let the compact and convex $K\subset \mathcal{M}^n$ with non-empty interior satisfy that
$\tau_H(K)$ is convex for any $(n-1)$-dimensional subspace $H$ of $\mathcal{M}^n$.

First we prove that for any pair $x,y\in \partial_{\mathcal{M}^n}K$, $x\neq y$, of strongly regular points,
writing $H$ to denote the perpendicular bisector $(n-1)$-subspace of $[x,y]_{\mathcal{M}^n}$, we have
\begin{equation}
\label{NK(x)}
N_K(x)=\sigma_H(N_K(y)).
\end{equation} 
Let  $x\in H^+$. 
To prove (\ref{NK(x)}), we verify that
\begin{equation}
\label{xonboundary}
x\in\partial_{\mathcal{M}^n} \tau_H(K)
\end{equation} 
by observing $[x,y]_{\mathcal{M}^n}\subset\tau_H(K)$ and distinguishing two cases.

If $[x,y]_{\mathcal{M}^n}\subset\partial_{\mathcal{M}^n} K$, then
there exists a supporting $(n-1)$-dimensional subspace $\Pi$ to $K$ containing $[x,y]_{\mathcal{M}^n}$
according to Lemma~\ref{secant}.
It follows that $\Pi$ is a  supporting $(n-1)$-dimensional subspace to $\sigma_H K$, and in turn
to $\tau_H (K)$, proving (\ref{xonboundary}) in this case.

On the other hand, if $[x,y]_{\mathcal{M}^n}\cap {\rm int}_{\mathcal{M}^n} K\neq \emptyset$, then let $\ell\subset \mathcal{M}^n$ be the one-dimensional subspace of $x$ and $y$, and hence
$\ell\cap K=[x,y]_{\mathcal{M}^n}$ by Lemma~\ref{secant}. As $x\in \partial_{\mathcal{M}^n} \sigma_H(K)$,
we deduce from Lemma~\ref{two-point-prop} (i) and (ii) that
$\ell\cap \tau_H(K) =[x,y]_{\mathcal{M}^n}$, and in turn (\ref{xonboundary}) finally follows.

Since $x\in\partial_{\mathcal{M}^n} \tau_H(K)$, there exists a supporting $(n-1)$-dimensional subspace
$\Xi$ to $\tau_H(K)$ at $x$, and let $\tau_H(K)\subset \Xi^+$. We deduce from Lemma~\ref{two-point-prop} (i) that $\sigma_H(K),K\subset \Xi^+$, thus $\sigma_H(\Xi)$ is a supporting $(n-1)$-dimensional subspace
to $K$ at $y$ with $K\subset\sigma_H(\Xi^+)$. As $x$ and $y$ are strongly regular points, we conclude
(\ref{NK(x)}) by Lemma~\ref{secant}.

In turn, we deduce from (\ref{NK(x)}) that
for any pair $x,y\in \partial_{\mathcal{M}^n}K$, $x\neq y$, of strongly regular points, there exists
$\lambda(x,y)\in\R$ such that
\begin{equation}
\label{lambda-def}
N_K(x)-N_K(y)=\lambda(x,y)(x-y)\in\R^{n+1}.
\end{equation}
Obviously, $\lambda(x,y)=\lambda(y,x)$.

We fix a strongly regular point $x_0\in  \partial_{\mathcal{M}^n}K$. We claim that if $x,y\in \partial_{\mathcal{M}^n}K$ are strongly regular points different from $x_0$, then
\begin{equation}
\label{x0xy}
\lambda(x_0,x)=\lambda(x_0,y).
\end{equation}
We distiguishing two cases. If  $x_0,x,y$ are not contained in a one-dimensional subspace, then  (\ref{lambda-def})
yields  that
\begin{eqnarray*}
N_K(x_0)-N_K(x)&=&\lambda(x_0,x)(x_0-x) \\
N_K(x)-N_K(y)&=&\lambda(x,y)(x-y) \\
N_K(y)-N_K(x_0)&=&\lambda(y,x_0)(y-x_0).
\end{eqnarray*}
Adding up the three relations yields 
$\lambda(x_0,x)(x_0-x)+\lambda(x,y)(x-y)+\lambda(y,x_0)(y-x_0)=o$.
Since any two of $x_0-x,x-y,y-x_0$ are independent in $\R^{n+1}$, we have
$\lambda(x_0,x)=\lambda(x,y)=\lambda(y,x_0)$, and hence (\ref{x0xy}) holds in this case.

On the other hand, if $x_0,x,y$ are contained in a one-dimensional subspace, then
as strongly regular points are dense in $\partial_{\mathcal{M}^n}K$, there exists a strongly regular point $z\in \partial_{\mathcal{M}^n}K$ not contained in the one-dimensional subspace passing through $x_0,x,y$. Applying the previous case
first to the triple $x_0,x,z$, then to the triple $x_0,y,z$, 
it follows that $\lambda(x_0,x)=\lambda(x_0,z)=\lambda(x_0,y)$, proving (\ref{x0xy}).

According to (\ref{x0xy}), there exists a common value $\lambda$ of $\lambda(x_0,x)$
 for all strongly regular points $x\in \partial_{\mathcal{M}^n}K$. Setting
$p=N_K(x_0)-\lambda x_0\in\R^{n+1}$, we deduce from (\ref{lambda-def}) that
\begin{equation}
\label{lambda-regpoints}
N_K(x)=p+\lambda  x \mbox{ \ for strongly regular point $x\in \partial_{\mathcal{M}^n}K$}.
\end{equation}

The rest of the argument is split between the hyperbolic- and spherical case.\\

\noindent{\bf Case 1} $\mathcal{M}^n=H^n$\\
We claim that
\begin{equation}
\label{pnotzeroHn}
p\neq o.
\end{equation}
Otherwise $N_K(x_0)=\lambda  x_0$ by (\ref{lambda-regpoints}), and hence
$$
-1=\mathcal{B}(N_K(x),N_K(x))=\lambda^2\mathcal{B}(x,x)=\lambda^2,
$$
what is absurd, verifying (\ref{pnotzeroHn}).

If $x\in \partial_{H^n}K$ is a strongly regular point, then $N_K(x)\in T_x$, thus (\ref{lambda-regpoints}) yields that 
$0=\mathcal{B}(N_K(x),x)=\mathcal{B}(p+\lambda  x,x)$, and hence 
$\mathcal{B}(p,x)=-\lambda$. As $p\neq o$, we deduce that each strongly regular point
of $\partial_{H^n}K$ is contained in the standard hypersurface
$\{x\in H^n:\,\mathcal{B}(p,x)=-\lambda\}$
(see the list before Lemma~\ref{hypersurf-component}), therefore
 Corollary~\ref{hypersurf-solid} yields that $K$ is a ball.\\

\noindent{\bf Case 2} $\mathcal{M}^n=S^n$\\
We again  claim that
\begin{equation}
\label{pnotzeroSn}
p\neq o.
\end{equation}
Otherwise $N_K(x_0)=\lambda  x_0$ by (\ref{lambda-regpoints}), and hence
$$
1=\langle N_K(x_0),N_K(x_0)\rangle=\langle N_K(x_0),\lambda  x_0\rangle=0,
$$
what is absurd, verifying (\ref{pnotzeroSn}).

If $x\in \partial_{S^n}K$ is a strongly regular point, then $N_K(x)\in T_x$, thus (\ref{lambda-regpoints}) yields that 
$0=\langle N_K(x),x\rangle=\langle p+\lambda  x,x\rangle$, and hence 
$\langle p,x\rangle=-\lambda$. As $p\neq o$, we deduce that each strongly regular point
of $\partial_{S^n}K$ is contained in the boundary
$\{x\in S^n:\,\langle p,x\rangle=-\lambda\}$ of a fixed spherical ball
(see the list before Lemma~\ref{hypersurf-component}), therefore
 Corollary~\ref{hypersurf-solid} yields that $K$ is a spherical ball.
\endproof

\noindent{\bf Remark: } In $\R^n$, a similar argument works, but some changes have to be instituted.
Instead of (\ref{NK(x)}), we have $N_K(x)=\sigma_{\widetilde{H}}(N_K(y))$ where
$\widetilde{H}$ is the linear $(n-1)$-plane in $\R^n$ parallel to $H$, therefore 
(\ref{lambda-def}) still holds. Similarly as above, (\ref{lambda-def}) leads to (\ref{lambda-regpoints}).
For the final part of the argument, we have $\lambda\neq 0$, because otherwise
(\ref{lambda-regpoints}) yields that $N_K(x)=p$ for any strongly regular point $x\in \partial_{R^n}K$, contradicting that the strongly regular 
points are dense in $x\in \partial_{R^n}K$. As (\ref{lambda-regpoints}) implies that
if  $x\in \partial_{R^n}K$ is a strongly regular point, then 
$$
\lambda^{-2}=\langle \lambda^{-1}N_K(x),\lambda^{-1}N_K(x)\rangle=
\langle \lambda^{-1}p+ x, \lambda^{-1}p+ x\rangle,
$$
and hence each strongly regular point
of $\partial_{\R^n}K$ is contained in the boundary
$\{x\in \R^n:\,\langle x+\lambda^{-1}p, x+\lambda^{-1}p\rangle=\lambda^{-2}\}$ of a fixed ball.
Therefore
 Corollary~\ref{hypersurf-solid} yields that $K$ is a ball.\\

\noindent{\bf Proof of Theorem~\ref{spherical} when $D\leq\frac{\pi}2$ and
of Theorem~\ref{hyperbolic} }
Let $\mathcal{M}^n=H^n$, or let $\mathcal{M}^n=S^n$ and $D\leq\frac{\pi}2$.
Theorem~\ref{maximal-sets} (i) yields the existence of $D$-maximal sets in $\mathcal{M}^n$,
and any $D$-maximal set is convex according to Corollary~\ref{maximal-convex}.

Let $C$ be any  $D$-maximal set, and hence $C$ is convex. We deduce from 
from Lemma~\ref{two-point-prop} that  for any $(n-1)$-dimensional subspace $H$ of
$\mathcal{M}^n$, $\tau_H(C)$
is $D$-maximal, thus convex, therefore Theorem~\ref{two-point-convex} yields that $C$ is a ball. The maximality of $V(C)$ implies that the radius of $C$ is $D/2$, proving Theorem~\ref{spherical} when $D\leq \pi/2$ and also
Theorem~\ref{hyperbolic}. $\Box$\\

\noindent{\bf Remark: } In $\R^n$, (\ref{eucklidean}) can proved in a similar way.

\section{Proof of Theorem~\ref{spherical} if $\frac{\pi}2<D<\pi$}
\label{secnonconvexproof}

Let  $\frac{\pi}2<D<\pi$. We frequently drop the index ``$S^n$" in the formulas, say we simply write $B(x,r)$, $V(X)$ and
$\partial X$ for $x\in S^n$, $X\subset S^n$ and $r\in(0,\pi)$.

Let us recall from Section~\ref{secDmaximal} that $C\subset S^n$ is $D$-maximal if
\begin{itemize}
\item $C$ is compact;
\item ${\rm diam}_{S^n}\,C\leq D$;
\item $V(C)=\sup\{V(X):\, X\subset S^n\mbox{ compact and }
{\rm diam}_{S^n} X\leq D \}$.
\end{itemize}
According to Theorem~\ref{maximal-sets} (i), there exist $D$-maximal sets in $S^n$.
We deduce from Lemma~\ref{two-point-prop} that if $C\subset S^n$ is $D$-maximal, then
\begin{equation}
\label{tauCmax}
\tau_H(C)
\mbox{ \ is $D$-maximal for any $(n-1)$-dimensional subspace $H$ of
$S^n$.}
\end{equation}

As $D>\frac{\pi}2$, it is {\it a priori} not clear whether a $D$-maximal set is convex. However,  
Theorem~\ref{maximal-sets} (ii) implies that for
 any $D$-maximal set $C$ in $S^n$ and $z\in\partial C$,
there exists $y\in C$ such that
\begin{equation}
\label{Dmaxinball}
C\subset B(y,D)\mbox{ \  and  $z\in \partial B(y,D)$.}
\end{equation}

\begin{lemma}
\label{solidmax}
If $\frac{\pi}2<D<\pi$ and $C$ is a $D$-maximal set in $S^n$, then there exists a solid $D$-maximal set 
$C_0\subset C$.
\end{lemma}
\proof Let $X\subset C$ be the set of density points; namely,
$$
X=\left\{z\in C: \lim_{r\to 0^+}\frac{V(C\cap B(z,r))}{V(B(z,r))}=1\right\},
$$
thus readily ${\rm int}\, C\subset X$. It follows from Lebesgue's Density Theorem that
$$
V(X)=V(C).
$$
In addition, if $z\in \partial C$, then (\ref{Dmaxinball}) yields that $z\not\in X$, and hence
$X={\rm int}\, C$. We conclude that $C_0$ can be taken as the closure of $X$.
\endproof

The idea of the proof of
Theorem~\ref{spherical} when $\frac{\pi}2<D<\pi$ is similar to the case $D\leq \frac{\pi}2$; more precisely, 
the idea is to prove that 
if $\frac{\pi}2<D<\pi$ and $C\subset S^n$ is a solid $D$-maximal set, then there exist $p\in\R^{n+1}$ and $\lambda\in\R$ such that 
\begin{equation}
\label{idea}
p+\lambda x\in\{N_C(x),-N_C(x)\} \mbox{ \ for any strongly regular point $x\in\partial C$.}
\end{equation}
Thus from now on, the main goal is to understand properties of solid $D$-maximal sets. We use two-point symmetrization again. The difficulty in the $D> \frac{\pi}2$ case is that 
if $x,y\in \partial C$, $x\neq y$, are strongly regular points of a solid $D$-maximal set $C$ in $S^n$, 
and $H$ is the perpendicular bisector of $[x,y]_{S^n}$, then
{\it a priori} $x$ may lie in ${\rm int}\,\tau_HC$.

\begin{lemma}
\label{pmNxy}
Let $\frac{\pi}2<D<\pi$, and let $C$ be a solid $D$-maximal set in $S^n$. 
If $x,y\in \partial C$, $x\neq y$, are strongly regular points, then
there exist
$\lambda(x,y)\in\R$ and $\eta(x,y)\in\{-1,1\}$ such that
\begin{equation}
\label{lambda-defSn}
N_K(x)-\eta(x,y)N_K(y)=\lambda(x,y)(x-y)\in\R^{n+1}
\end{equation}
where $\eta(y,x)=\eta(x,y)$ and $\lambda(y,x)=\eta(x,y)\lambda(x,y)$.
\end{lemma}
\proof It is equivalent to prove that
\begin{equation}
\label{NK(x)S^n}
N_C(x)=\pm\sigma_H(N_C(y))
\end{equation} 
where $H$ is the perpendicular bisector $(n-1)$-dimensional subspace of $[x,y]_{S^n}$. In turn,
 (\ref{NK(x)S^n}) is equivalent proving that if we assume
\begin{equation}
\label{NK(x)S^nno}
N_C(x)\neq -\sigma_H(N_C(y)),
\end{equation} 
then we have
\begin{equation}
\label{NK(x)Sn}
N_C(x)=\sigma_H(N_C(y)).
\end{equation} 
Let  $x\in H^+$. We deduce from (\ref{Dmaxinball}) and from the fact that $x$ and $y$ are strongly regular boundary points
that 
there exist $r>0$ and $x_0,x_1,y_0,y_1\in C$ such that
\begin{eqnarray*}
B(x_0,r)&\subset C\subset&B(x_1,D) \mbox{ \ where 
$x\in \partial B_{S^n}(x_0,r)\cap \partial B_{S^n}(x_1,D)$ } \\
B(y_0,r)&\subset C\subset&B(y_1,D) \mbox{ \ where 
$y\in \partial B_{S^n}(y_0,r)\cap \partial B_{S^n}(y_1,D)$} \\
B(x_0,r)\cap H&=&B(y_0,r)\cap H=\emptyset.
\end{eqnarray*}
Readily, we have $N_C(x)=N_{B(x_0,r)}(x)=N_{B(x_1,D)}(x)$ and
$N_C(y)=N_{B(y_0,r)}(y)=N_{B(y_1,D)}(y)$.

It follows from (\ref{tauCmax}) that $\tau_HC$ is $D$-maximal as well.
As $x=\sigma_Hy\in \tau_HC$ and
$$
\tau_HC\subset B(x_1,D)\cup \sigma_H B(y_1,D)
$$ 
by Lemma~\ref{two-point-prop}, and
(\ref{NK(x)S^nno}) yields that $x$ is a boundary point of $B(x_1,D)\cup \sigma_H B(y_1,D)$,
we deduce that 
$x\in\partial \tau_HC$. We deduce from 
(\ref{Dmaxinball}) that
$\tau_HC\subset B(z,D)$ and  $x\in \partial B(z,D)$ for some $z\in\tau_HC$.
Now $x\in B(x_0,r)\cap \sigma_HB(y_0,r)\cap \partial B(z,D)$ and
$$
B(x_0,r)\subset C\subset B(z,D)\mbox{ \ and \ }\sigma_HB(y_0,r)\subset C\subset B(z,D)
$$
by Lemma~\ref{two-point-prop}, therefore $B(x_0,r)=\sigma_HB(y_0,r)$. We conclude 
(\ref{NK(x)Sn}), and in turn (\ref{NK(x)S^n}), proving Lemma~\ref{pmNxy}.
\endproof

As a first step to prove (\ref{idea}), we consider certain specific triples of strongly regular points. 

\begin{lemma}
\label{triple}
Let $\frac{\pi}2<D<\pi$, and let $C$ be a solid $D$-maximal set in $S^n$. 
If $x_0,x,y\in \partial C$ are strongly regular points not contained in a one-dimensional subspace of $S^n$
such that $\eta(x_0,y)=1$ and $\eta(x_0,x)=\eta(y,x)$, then for $\lambda=\lambda(y,x_0)=\lambda(x_0,y)$, we have
$\lambda(x_0,x)=\lambda$ and
$$
\eta(x_0,x) N_C(x)-\lambda x=N_C(x_0)-\lambda x_0.
$$
\end{lemma}
\proof  For $\eta(x_0,x)=\eta(y,x)=\eta$, Lemma~\ref{pmNxy} implies that 
\begin{eqnarray*}
N_C(x)-\eta N_C(y)&=&\lambda(x,y)(x-y) \\
N_C(x_0)-\eta N_C(x)&=&\lambda(x_0,x)(x_0-x) \\
N_C(y)-N_C(x_0)&=&\lambda(y,x_0)(y-x_0).
\end{eqnarray*}
Since $\eta\in\{-1,1\}$, we replace the first equality by $\eta N_C(x)- N_C(y)=\eta \lambda(x,y)(x-y)$
and add up these three relations, and hence we obtain
$$
\eta \lambda(x,y)(x-y)+\lambda(x_0,x)(x_0-x)+\lambda(y,x_0)(y-x_0)=o.
$$
Since any two of $x_0-x,x-y,y-x_0$ are independent in $\R^{n+1}$, we have
$\eta \lambda(x,y)=\lambda(x_0,x)=\lambda(y,x_0)$ where $\lambda(y,x_0)=\lambda$.
We deduce that
$N_C(x_0)-\eta N_C(x)=\lambda (x_0-x)$, proving Lemma~\ref{triple}.
\endproof

Let us show that we have the setup in Lemma~\ref{triple} if $x$ and $y$ are close enough strongly regular boundary points.

\begin{lemma}
\label{close-strongly regular}
Let $\frac{\pi}2<D<\pi$, and let $C$ be a solid $D$-maximal set in $S^n$, and let $x_0\in \partial C$ be a strongly regular point. 
\begin{description}
\item[(i)] If $\lim_{m\to \infty}y_m=x_0$ for $y_m\in\partial C$ strongly regular boundary points, then
$$
\lim_{m\to \infty}N_C(y_m)=N_C(x_0)\mbox{ \ and \ }
\lim_{m\to \infty}\left\langle N_C(x_0),\frac{y_m-x_0}{d_{\R^{n+1}}(y_m,x_0)}\right\rangle=0;
$$ 
\item[(ii)] For any strongly regular point $z\in \partial C$, $z\neq x_0$, there exists $\varrho\in(0,\frac{\pi}2)$ depending on $x_0,z$ and $C$ such that
$\eta(x_0,y)=1$ and $\eta(x_0,z)=\eta(y,z)$ if
$y\in B(x_0,\varrho)\cap\partial C$ is a strongly regular point with $y\neq x_0$.
\end{description}
\end{lemma}
\proof There exist some $r\in(0,\frac{\pi}2)$ and $z_0,z_1\in C$ such that 
$B(z_0,r)\subset C\subset B(z_1,D)$ and $x_0\in \partial B(z_0,r)\cap\partial B(z_1,D)$
according to Theorem~\ref{maximal-sets} (ii).

For (i), there exist some $w_m\in C$  with $C\subset B(w_m,D)$ and $y_m\in\partial B(w_m,D)$, and we may assume that $B(w_m,D)$ tends to some $B(w,D)$ for $w\in C$. As $B(z_0,r)\subset B(w,D)$ and
$x_0\in \partial B(z_0,r)\cap\partial B(w,D)$, we have $w=z_1$. This yields that
$\lim_{m\to \infty}N_C(y_m)=N_C(x_0)$. In addition, 
$\lim_{m\to \infty}\left\langle N_C(x_0),\frac{y_m-x_0}{d_{\R^{n+1}}(y_m,x_0)}\right\rangle=0$ follows from
$\lim_{m\to \infty}y_m=x_0$ and $y_m\in B(z_1,D)\backslash{\rm int}\,B(z_0,r)$.

We prove (ii) by contradiction, and we assume that there exists a sequence of
strongly regular boundary points $y_m\in\partial C$ such that $\lim_{m\to \infty}y_m=x_0$, and
\begin{description}
\item[(a)] either $N_C(x_0)+N_C(y_m)=\lambda(x_0,y_m)(x_0-y_m)$ for each $m$,
\item[(b)] or $N_C(z)+\eta(z,x_0)N_C(y_m)=\lambda(z,y_m)(z-y_m)$ for each $m$.
\end{description}
If (a) holds, then readily $|\lambda(x_0,y_m)|\leq 2/d_{\R^{n+1}}(x_0,y_m)$, thus
(i) yields that
\begin{eqnarray*}
4&=&\lim_{m\to \infty}\langle N_C(x_0)+N_C(y_m),N_C(x_0)+N_C(y_m)\rangle=
\lim_{m\to \infty}\left|\langle N_C(x_0)+N_C(y_m),\lambda(x_0,y_m)(x_0-y_m)\rangle\right|\\
&\leq & 
2\lim_{m\to \infty}\left|\left\langle N_C(x_0)+N_C(y_m),\frac{x_0-y_m}{d_{\R^{n+1}}(y_m,x_0)}\right\rangle\right|
=0,
\end{eqnarray*}
what is absurd. 

If (b) holds, then as $z-x_0$ and $N_C(x_0)\in T_{x_0}$ are independent, there exists some
$u\in S^{n-1}$ such that $\langle u, z-x_0\rangle=0$ and $\langle u, N_C(x_0)\rangle >0$. We also observe that
 if $m$ is large, then
$|\lambda(z,y_m)|\leq 3/d_{\R^{n+1}}(z,x_0)$.

We deduce from $\lim_{m\to \infty}N_C(y_m)=N_C(x_0)$ that
\begin{eqnarray*}
2\langle u, N_C(x_0)\rangle&=&\lim_{m\to \infty}\left|\langle u,[N_C(z)+\eta(z,x_0)N_C(x_0)]- [N_C(z)-\eta(z,x_0)N_C(x_0)]\rangle\right|\\
&=&
\lim_{m\to \infty}\left|\langle u,\lambda(z,y_m)(z-y_m)-\lambda(z,x_0)(z-x_0)\rangle\right|\\
&= & 
\lim_{m\to \infty}|\lambda(z,y_m)|\cdot \left|\langle u,(z-y_m)\rangle\right|\leq
\frac{3}{d_{\R^{n+1}}(z,x_0)}\cdot\lim_{m\to \infty}\left|\langle u,(z-y_m)\rangle\right|
=0,
\end{eqnarray*}
what is again a contradiction, proving (ii). 
\endproof

Now we choose the right  strongly regular ``base point" $x_0$.

\begin{lemma}
\label{base-point}
Let $\frac{\pi}2<D<\pi$, and let $C$ be a solid $D$-maximal set in $S^n$. There exists
 a strongly regular point $x_0\in \partial C$ such that for any 
one-dimensional subspace $\ell$ of $S^n$ passing through $x_0$ and any $\varepsilon\in(0,\frac{\pi}2)$,
one finds a strongly regular point $y\in \partial C\cap(B(x_0,\varepsilon)\backslash \ell)$.
\end{lemma}
\proof First let $n\geq 3$, let $x_0\in \partial C$ be any strongly regular point, and let $\varepsilon\in(0,\frac{\pi}2)$.
As $C$ is solid, there exist 
$$
z\in {\rm int}B(x_0,\varepsilon)\cap [C\backslash \ell]\mbox{ \ and \ } 
y\in {\rm int}B(x_0,\varepsilon)\backslash [C\cup \ell].
$$ 
Since ${\rm int}B(x_0,\varepsilon)\backslash \ell$ is connected (this is the point where we use that $n\geq 3$), connecting $y$ and $z$ by
a continuous curve in ${\rm int}B(x_0,\varepsilon)\backslash \ell$ implies that there exists a
$w\in \partial\, C\cap [{\rm int}B(x_0,\varepsilon)\backslash\ell]$. It follows from the fact
that strongly regular points in $\partial C$ are dense that there exists
a strongly regular point $\tilde{w}\in \partial\, C\cap [{\rm int}B(x_0,\varepsilon)\backslash\ell]$.

Next let $n=2$. In this case, the argument is indirect. Since strongly regular points are dense on the boundary,
we suppose that for any strongly regular point $w\in\partial C$, there exist a one-dimensional subspace $\ell_w$
pasing through $w$
and an $r_w\in(0,\frac{\pi}2)$ such that
\begin{equation}
\label{lineindisc0}
\partial C\cap B(w,r_w)\subset \ell_w,
\end{equation}
and seek a contradiction. As $C$ is solid, there exist 
$$
z\in {\rm int}B(w,r_w)\cap [C\backslash \ell_w]\mbox{ \ and \ } 
y\in {\rm int}B(w,r_w)\backslash [C\cup \ell_w].
$$ 
For any $x\in B(w,r_w)\cap \ell_w$, the piecewise linear path $[y,x]_{S^n}\cup [x,z]_{S^n}$
intersects $\partial C$, and the intersection can be only $x$ by (\ref{lineindisc0}). Therefore
\begin{equation}
\label{lineindisc}
\partial C\cap B(w,r_w)= \ell_w\cap B(w,r_w).
\end{equation}
It also follows that one component of ${\rm int}B(w,r_w)\backslash \ell_w$ is part of ${\rm int}C$ and
the other component is disjoint from $C$,
and hence any $x\in \ell_w\cap {\rm int}B(z,r_w)$ is a strongly regular boundary point with
\begin{equation}
\label{lineindiscnormal}
N_C(x)=N_C(w)\mbox{ \ for each $x\in \ell_w\cap {\rm int}B(z,r_w)$ and \ }
l_w=\{x\in S^n:\,\langle x,N_C(w)\rangle=0\}.
\end{equation}

Let us choose a strongly regular point $w\in\partial C$ and an other strongly regular point 
$v\in \partial C\backslash \ell_w$. It follows that $\ell_v\neq \ell_w$.
We deduce from Lemma~\ref{close-strongly regular} (ii) and (\ref{lineindiscnormal}) that
there exists a $y\in l_w\cap\partial C$ with $y\neq w$ and $\eta\in\{-1,1\}$ such that
$$
N_C(y)=N_C(w),\mbox{ \ }\eta(y,w)=1\mbox{ \ and \ }
\eta=\eta(v,w)=\eta(v,y).
$$
On the one hand, we deduce from
from $N_C(y)=N_C(w)$ that $\lambda(w,y)=\lambda(y,w)=0$ in Lemma~\ref{pmNxy}, therefore Lemma~\ref{triple} yields that $\eta N_C(v)=N_C(w)$. 
However  $\ell_v\neq \ell_w$ and (\ref{lineindiscnormal}) imply
$N_C(v)\neq \pm N_C(w)$ in $\R^{n+1}$, thus we have arrived at a contradiction,
verifying Lemma~\ref{base-point} also if $n=2$.
\endproof

\noindent{\bf Proof of Theorem~\ref{spherical} when $\mathbf{\frac{\pi}2<D<\pi}$: }
Let $C$ be any  $D$-maximal set, and let $C_0\subset C$ be the solid $D$-maximal set provided by
Lemma~\ref{solidmax}. According to Lemma~\ref{base-point}, there exists
 a strongly regular point $x_0\in \partial C_0$ such that for any 
one-dimensional subspace $\ell$ passing through $x_0$ and any $\varepsilon\in(0,\frac{\pi}2)$,
one finds a strongly regular point $y\in \partial C\cap(B(x_0,\varepsilon)\backslash \ell)$.

We deduce from Lemma~\ref{close-strongly regular} (ii)
that there exists $\tilde{r}\in(0,\frac{\pi}2)$ depending on $x_0$ and $C$
 such that
$\eta(x_0,y)=1$ for any strongly regular point
$y\in B(x_0,\tilde{r})\cap\partial C$ with $y\neq x_0$. We claim that there exist
$r\in(0,\tilde{r})$ and
$\lambda\in\R$ such that
\begin{equation}
\label{yclosetox0}
\lambda(x_0,y)=\lambda(y,x_0)=\lambda\mbox{ \ and \ } 
\eta(x_0,y)=1\mbox{ holds for any 
$y\in B(x_0,r)\cap\partial C$ with $y\neq x_0$.}
\end{equation}
Readily, $\eta(x_0,y)=1$ by $r<\tilde{r}$. To have the right value of $\lambda(x_0,y)$,
first we fix a strongly regular point
$y_0\in B(x_0,\tilde{r})\cap\partial C$ with $y_0\neq x_0$, and write $\ell_0$ to denote
the one-dimensional subspace spanned by $x_0$ and $y_0$. Set
$$
\lambda=\lambda(x_0,y_0).
$$
Lemma~\ref{triple} and Lemma~\ref{close-strongly regular} (ii) applied to the triple $x_0,y,y_0$ imply the existence of an $r_0\in(0,\tilde{r})$
such that $\lambda(x_0,y)=\lambda$ holds for any strongly regular point
$y\in \partial C\cap(B(x_0,r_0)\backslash \ell_0)$. 
We fix such a strongly regular point 
$y_1\in \partial C\cap(B(x_0,r_0)\backslash \ell_0)$, and write $\ell_1$ to denote
the one-dimensional subspace spanned by $x_0$ and $y_1$. In particular,
$\ell_0\neq \ell_1$.

Finally, 
applying  Lemma~\ref{triple} and Lemma~\ref{close-strongly regular} (ii)
to the triple $x_0,y,y_1$, there exists
 $r\in(0,r_0)$
such that $\lambda(x_0,y)=\lambda$ holds for any strongly regular point
$y\in \partial C\cap(B(x_0,r)\backslash \ell_1)$. Since 
either $y\not\in\ell_0$ or $y\not\in\ell_1$ hold for any point
$y\in \partial C\cap(B(x_0,r)\backslash x_0)$, we conlcude (\ref{yclosetox0}).

Our next goal is to verify that $p=N_C(x_0)-\lambda x_0$ satisfies that
if $x\in\partial C$ is a strongly regular point, then
\begin{equation}
\label{pNCx}
p+\lambda x=\eta(x_0,x)N_C(x)
\end{equation}
where we set $\eta(x_0,x_0)=1$. Writing $\ell$ to denote the one-dimensional subspace
passing through $x_0$ and $x$, the choice of $x_0$ and Lemma~\ref{close-strongly regular} (ii)
yield the existence of a strongly regular point $y\in \partial C\cap(B(x_0,r)\backslash \ell)$
such that $\eta(x_0,y)=1$ and $\eta(x_0,x)=\eta(y,x)$.
Therefore we conclude (\ref{pNCx}) by Lemma~\ref{triple}.

We again  claim that
\begin{equation}
\label{pnotzeroSn0}
p\neq o.
\end{equation}
Otherwise $N_K(x_0)=\lambda  x_0$ by (\ref{pNCx}), and hence
$$
1=\langle \eta(x_0,x)N_K(x_0),\eta(x_0,x)N_K(x_0)\rangle=\langle \eta(x_0,x)N_K(x_0),\lambda  x_0\rangle=0,
$$
what is absurd, verifying (\ref{pnotzeroSn0}).

If $x\in \partial C$ is a strongly regular point, then $N_K(x)\in T_x$, thus (\ref{pNCx}) yields that 
$0=\langle \eta(x_0,x)N_{C_0}(x),x\rangle=\langle p+\lambda  x,x\rangle$, and hence 
$\langle p,x\rangle=-\lambda$. As $p\neq o$, we deduce that each strongly regular point
of $\partial C_0$ is contained in the boundary
$\{x\in S^n:\,\langle p,x\rangle=-\lambda\}$ of a fixed spherical ball, therefore
 Corollary~\ref{hypersurf-solid} yields that $C_0$ is a spherical ball.
 As $C_0$ has maximal volume among sets of diameter at most $D$,
it follows that the radius of $C_0$ is $D/2$, say $C_0=B(z,D/2)$.

Finally, we show that  $C=B(z,D/2)$. To prove this, let $x\in S^n\backslash  B(z,D/2)$, and
let $\ell$ be the (or a) one-dimensional subspace of $S^n$ passing through $x$ and $x_0$. As $\ell$ intersects
$B(z,D/2)$ in an arc of length $D$, this arc contains a point $y$ with $d_{S^n}(x,y)>D$. Therefore
$x\not \in C$, completing the proof of Theorem~\ref{spherical} when $D> \pi/2$. $\Box$ \\

\noindent{\bf Acknowledgement: } We are grateful for the careful referee whose remarks have substantially improved the paper. 
K\'aroly J. B\"or\"oczky is supported by NKFIH projects ANN 121649, K 109789,
 K 116451 and KH 129630.

\end{document}